\documentstyle[12pt]{article}
\newcommand{\const}{\mathop{\rm const}\limits}

\newcommand{\supp}{\mathop{\rm supp}\limits}

\newcommand{\Var}{\mathop{\rm Var}\limits}

\newcommand{\Law}{\mathop{\rm Law}\limits}

\newcommand{\Cov}{\mathop{\rm Cov}\limits}

\begin{document}

\begin{center}

\vspace{3mm}

{\bf Well Posedness  of the Problem of Estimation }\par

\vspace{4mm}

{\bf Fractional Derivative for a  Distribution Function.} \\

\vspace{4mm}

 $ {\bf E.Ostrovsky^a, \ \ L.Sirota^b } $ \\

\vspace{4mm}

$ ^a $ Corresponding Author. Department of Mathematics and computer science, Bar-Ilan University, 84105, Ramat Gan, Israel.\\

E-mail: eugostrovsky@list.ru\\

\vspace{3mm}

$ ^b $  Department of Mathematics and computer science. Bar-Ilan University,
84105, Ramat Gan, Israel.\\

E-mail: sirota3@bezeqint.net \\

\vspace{4mm}
                    {\sc Abstract.}\\

 \end{center}

 \vspace{3mm}

  We study the problem of nonparametric estimation of the fractional derivative of unknown distribution function
and of spectral function and show that these problems are well posed when the order of derivative is less than 0.5. \par
  We prove also the unbiaseness  and asymptotical normality of offered estimates with optimal speed of convergence.\par
   For the construction of the  confidence region  in some functional norm we establish the Central Limit Theorem in
 correspondent Lebesgue-Riesz space  for offered estimates, and deduce also the non-asymptotical deviation of our
 estimates in these spaces.  \par

 \vspace{4mm}

{\it Key words and phrases:} Fractional derivatives and integrals of a Rieman-Liouville type, empirical and exact function of
distribution, reliability function, loss functional, indicator function, density, spectral function and density,
sample, estimate, confidence region, periodogram,
asymptotical normality, bias  and unbiased estimate,  Gaussian random process, Kolmogorov's theorem, Central Limit Theorem
in Banach space, Lebesgue-Riesz and Grand Lebesgue spaces (GLS), measurable set, random variable (r.v.)
and random process (r.p.),  measurable function. \par

\vspace{3mm}

{\it Mathematics Subject Classification (2000):} primary 60G17; \ secondary 60E07; 60G70.\\

\vspace{4mm}

\section{Notations. Statement of problem.}

\vspace{4mm}

 "Fractional derivatives have been around for centuries  but recently they have
found new applications in physics, hydrology and finance", see  \cite{Meerschaert1}. \par
 Another applications: in the theory of Differential Equations are described in \cite{Miller1};
in statistics see in \cite{Adler1}, \cite{Bapna1}, see also \cite{Golubev1}, \cite{Enikeeva1}; in the theory of
integral equations etc. see in  the classical monograph \cite{Samko1}.  \par

\vspace{3mm}

 {\bf  We consider here  the problem of the nonparametric  estimation of the
fractional derivative for a  distribution function based on the sample of a "great" volume,
and analogously estimation of the fractional derivative of the spectral function of Gaussian stationary sequence. } \par

\vspace{3mm}

 We will prove that if the order of the fractional derivative $  \alpha $ is less than 1/2,  then these problems
are well posed. In particular, the speed of convergence  of offered unbiased estimate is $ 1/\sqrt{n}, $ as in the case of estimation
of ordinary distribution function $  F(x); $ they are asymptotical normal still in some rearrangement invariant norm. \par
 Our results improve ones in the articles  \cite{Borla1}, \cite{Golubev1},  \cite{Enikeeva1}, \cite{Matsui1} etc.,
 but does not contradict to the known results.\\

\vspace{3mm}

 More detail description. Let $  \xi_1, \xi_2, \ldots, \xi_n $ be a sample of a volume $ n, $ i.e. independent identical distributed
numerical random variable with common distribution function $  F = F(x). $ In what follows we restrict ourselves by consideration
of the following class $  K  $ of distributions:

$$
\forall F \in K \ \Rightarrow  \exists (a,b) \in R^2, \ 0 \le a < b \le \infty, \ F(a+0) = 0, F(b-0) = 1,
$$
such that on the interval $ (a,b) $ the function  $ F(x) $ is continuous and strictly increasing. \par

\vspace{3mm}

 We can and will suppose further without loss of generality $  a = 0. $ \par

  Let $ \alpha = \const \in (0,1); $ and let $ g = g(x), \  x \in R $  be measurable numerical function. The fractional derivative
 of a Rieman-Liouville type of order $  \alpha: \  D^{\alpha}[g](x)  = g^{(\alpha)}(x) $  is defined as follows:
 $ \Gamma(1-\alpha)  g^{(\alpha)}(x) =  $

$$
\Gamma(1-\alpha) \ D^{\alpha}[g](x)  =
\Gamma(1-\alpha) \ D^{\alpha}_x[g](x) \stackrel{def}{= } \frac{d}{dx} \int_0^x \frac{g(t) \ dt}{(x-t)^{\alpha}}. \eqno(1.1)
$$
see, e.g. the classical monograph of S.G.Samko, A.A.Kilbas and O.I.Marichev \cite{Samko1}, pp. 33-38; see also \cite{Miller1}.\par

 Hereafter $ \Gamma(\cdot) $ denotes the ordinary $ \Gamma \ $ function. \par

 We agree to take $ D^{\alpha}[g](x_0) = 0,  $ if at the point $  x_0 $ the expression $  D^{\alpha}[g](x_0)  $ does not exists. \par

Notice that the operator of the fractional derivative is non - local, if $ \alpha $ is not integer non-negative number.\par

 Recall also that the fractional  integral $  I^{(\alpha)}[\phi](x) = I^{\alpha}[\phi](x)  $ of a Rieman - Liouville type of
an order $ \alpha, 0 < \alpha < 1 $ is defined as follows:

$$
I^{(\alpha)}[\phi](x) \stackrel{def}{=} \frac{1}{\Gamma(\alpha)} \cdot \int_0^x \frac{\phi(t) \ dt}{(x-t)^{1 - \alpha}}, \ x,t > 0. \eqno(1.1a)
$$
 It is known  (theorem of Abel, see \cite{Samko1}, chapter 2, section 2.1)
 that the operator $ I^{(\alpha)}[\cdot]  $ is inverse to the fractional derivative operator $ D^{(\alpha)}[\cdot],  $
at least  in the class of absolutely continuous functions. \par

  Note that for the considered further functions this fractional derivative there exists almost everywhere. \par

\vspace{3mm}

 Let us consider the following important example. Define the function

$$
g_h(x) = I(h < x), \ x > 0, \ h = \const > 0.
$$

 We conclude after simple calculations taking into account our agreement

$$
g_h^{(\alpha)}(x) =  \frac{1}{\Gamma(1-\alpha)} \cdot I(h < x) \cdot (x - h)^{-\alpha}, \ \alpha = \const \in (0,1).
$$

 Let us calculate for the verification the fractional integral $  I^{\alpha}  $ of order $ \alpha $ from the function
$ x \to g_h^{(\alpha)}(\cdot). $ We have

$$
\Gamma(\alpha) \Gamma(1 - \alpha)  I^{\alpha}  \left[ g_h^{(\alpha)} \right] (x) = 0, \ x \le h,
$$
and in the case $  x > h $

$$
\Gamma(\alpha) \ \Gamma(1 - \alpha) \  I^{\alpha}  \left[ g_h^{(\alpha)} \right] (x) =
\int_0^x \frac{I(t > h) \ (t - h)^{-\alpha} \ dt}{ (x - t)^{1 - \alpha}} =
$$

$$
\int_h^x (t-h)^{-\alpha} \ (x - t)^{\alpha - 1} \ dt.
$$
 We make the substitution

 $$
 t = h + y(x-h); \hspace{6mm} dt = (x-h) dy:
 $$

$$
\Gamma(\alpha) \ \Gamma(1 - \alpha) \  I^{\alpha}  \left[ g_h^{(\alpha)} \right] (x) =
\int_0^1 y^{ - \alpha} \ (1 - y)^{\alpha - 1} \ dy =
$$

$$
B(1 - \alpha, \alpha) = \Gamma(\alpha) \ \Gamma(1 - \alpha) / \Gamma(1) =   \Gamma(\alpha) \ \Gamma(1 - \alpha),
$$
where $ B(\cdot, \cdot) $ denotes the usually Beta function. \par
 Thus,

$$
I^{\alpha}  \left[ g_h^{(\alpha)} \right] (x)  =  I(h < x)  =   g_h(x).
$$

 Note that since the function $ x \to  g_h(x) $ is not absolutely continuous, this result  can not be obtained from the
results of chapter 2 from the monograph \cite{Samko1}. \par

\vspace{5mm}

Further, we define as the capacity of a loss function  the following $ L_q(R, dF)  $ functional

$$
\tilde{W}_{n,q} [ D^{\alpha}[F](\cdot), \tilde{F}_{\alpha,n}(\cdot) ] \stackrel{def}{=} \sqrt{n} \times
\left[ {\bf E} \int_R | D^{\alpha}[F](x) - \tilde{F}_{\alpha,n}(x)|^q \ dF(x) \right]^{1/q}, \eqno(1.2)
$$
where $ q = \const \ge 1, \ \tilde{F}_{\alpha,n}(x) $  is arbitrary estimation of  $ D^{\alpha}[F](\cdot) $ based on our sample. \par
 But it is more convenient sometimes  to consider the equivalent problem of estimation of the fractional derivative of
 so-called "reliability" function $ D^{\alpha}[G](x), $ where

 $$
 G(x) = {\bf P} (\xi_i \ge x) = 1 - F(x), \eqno(1.3)
 $$
and to take $ G^{(\alpha)}(x) = D^{\alpha}G(x),  $

$$
W_{n,q} [ D^{\alpha}[G](\cdot), \tilde{G}_{\alpha,n}(\cdot) ] \stackrel{def}{=}  \sqrt{n} \times
\left[ {\bf E} \int_R | D^{\alpha}[G](x) - \tilde{G}_{\alpha,n}(x)|^q \ dF(x) \right]^{1/q}, \eqno(1.4)
$$
where in turn $ \tilde{G}_{\alpha,n}(x) $  is arbitrary estimation of  $ D^{\alpha}[G](\cdot) $ based on our sample. \par
 Note that

 $$
 D^{\alpha}[1](x) = \frac{x^{-\alpha}}{\Gamma(1 - \alpha) } \ne 0,
 $$
so that

$$
 D^{\alpha}[G](x) = \frac{x^{-\alpha}}{\Gamma(1 - \alpha) } -  D^{\alpha}[F](x) \ne -  D^{\alpha}[F](x). \eqno(1.4a)
$$

 For the practical using  the expression  (1.4) may be consistent approximate as $ n \to \infty $ as follows

 $$
W_q [ D^{\alpha}[G](\cdot), \tilde{G}_{\alpha,n}(\cdot) ] \approx \sqrt{n} \times
\left[ {\bf E} \int_R | D^{\alpha}[G](x) - \tilde{G}_{\alpha,n}(x)|^q \ dF_n(x) \right]^{1/q},
$$
 where $  F_n(x) $ is ordinary empirical function of distribution. \par
 We can define analogously the following estimate of the function $  G(x): $

 $$
 G_n(x) := n^{-1} \sum_{i=1}^n I(\xi_i \ge x), \eqno(1.5)
 $$
empirical reliability function. Here $ I(\xi_i \ge x) $ is the usually  indicator function:

$$
I(\xi_i \ge x) = 1 \ \Leftrightarrow \xi_i \ge x; \hspace{6mm}  I(\xi_i \ge x) = 0 \ \Leftrightarrow \xi_i < x.
$$

  Evidently, (Kolmogorov's theorem), the problem of distribution function estimation $ (\alpha = 0)  $ is well posed.
 V.D.Konakov in \cite{Konakov1}  proved in contradiction that  the problem of density estimation, i.e. when $ \alpha = 1, $ is ill posed. \par
  Roughly speaking, the result of V.D.Konakov  may be reformulated as follows. Certain problem of statistical estimation is well posed
iff there exists an estimate (more exactly, a sequence  of estimates) such that  the speed of convergence is equal (or less
than) $  1/\sqrt{n}. $ As a rule these estimations are asymptotically normal.\par

\vspace{4mm}

\section{Point estimate.}

\vspace{4mm}

{\bf 0.} {\it We suppose in what follows in this section  that $  x > 0 $ and } $ 0 < \alpha < 1/2, $
so that $ a = 0, \ 0 < b \le \infty. $  \par

\vspace{3mm}

{\bf 1.} Let us consider the following function

$$
x \to f_h(x) = I(x < h), \ h = \const > 0. \eqno(2.1)
$$
 It is easy to calculate that

$$
f_{\alpha,h}(x) := \Gamma(1 - \alpha) \ D^{\alpha}[f_h](x)  = x^{-\alpha} - (x - h)^{-\alpha} \cdot I(x > h), \ x > 0. \eqno(2.2)
$$

\vspace{3mm}

{\bf 2.}  It is reasonable to offer as a capacity  of the estimate $ G_{\alpha,n}(x) $ at the fixed point $  x, \ x > 0 $
of the fractional derivative $ D^{\alpha}[G](x)  $ the following statistic:

$$
\Gamma(1-\alpha) \cdot G_{\alpha,n}(x) \stackrel{def}{=} n^{-1} \sum_{i=1}^n f_{\alpha,\xi_i}(x). \eqno(2.3)
$$

 Denote also

$$
\Sigma^2_{\alpha}(x) := 2x^{-\alpha} \Gamma(1 - \alpha) G^{(\alpha)}(x) - \Gamma(1 - 2 \alpha) \ G^{(2 \alpha)}(x) -
\left(\Gamma(1 - \alpha) \ G^{(\alpha)}(x) \right)^2. \eqno(2.4)
$$

\vspace{3mm}

{\bf 3. \ Theorem 2.1.} \\
{\bf A.} Suppose that at the fixed positive point $ x $ the fractional derivative  $ G^{(\alpha)}(x) $ there exists.
Then the estimate $ G_{\alpha,n}(x) $ is not biased:

$$
 {\bf E} G_{\alpha,n}(x) = G^{(\alpha)}(x). \eqno(2.5)
$$

\vspace{3mm}

{\bf B.} Suppose in addition that at the fixed positive point $ x $
the fractional derivative  $ G^{(2 \alpha)}(x) $ there exists.  Then the estimate $ G_{\alpha,n}(x) $
is asymptotically as $  n \to \infty $ normal with the variance

$$
\Var \left[ G_{\alpha,n}(x) \right] = \frac{\Sigma^2_{\alpha}(x)}{n \ \Gamma^2(1 - \alpha)}: \eqno(2.6)
$$

\vspace{3mm}

$$
\Law \left\{ G_{\alpha,n}(x) \right\} \sim N \left(G^{(\alpha)}(x), \ \frac{\Sigma^2_{\alpha}(x)}{n \ \Gamma^2(1 - \alpha)} \right). \eqno(2.7)
$$

\vspace{3mm}

{\bf Proof.} It is sufficient to prove the equality (2.5) only for the value $  n = 1. $
 We deduce by direct computation  using the source definition (1.1)

$$
\Gamma(1 - \alpha) \cdot D^{\alpha}_x I(\xi \ge x)  \ = \ x^{-\alpha} - (x-\xi)^{-\alpha} \cdot I(\xi < x)=
$$

$$
f_{\alpha,\xi}(x), \hspace{6mm} \xi = \xi_1. \eqno(2.8)
$$
 It remains to take the expectation $ {\bf E} $ from both the sides of the relationship (2.8) to establish the unbiaseness. \par
 Let us calculate  the variance; we consider of course  the case $  n = 1. $

$$
\Var := \Var \left[ x^{-\alpha} - (x - \xi)^{-\alpha} I(\xi < x)   \right] =
$$
$$
\Var \left[ (x - \xi)^{-\alpha} I(\xi < x)   \right] = S_2 - S_1^2, \eqno(2.9)
$$
where

$$
S_1 = {\bf E} (x - \xi)^{-\alpha} \ I(\xi < x)
$$
and we know that

$$
S_1 = x^{-\alpha} - \Gamma(1 - \alpha) \ G^{(\alpha)}(x). \eqno(2.10)
$$

 Further, we will use the formula (2.24), section 2, pp. 35-37 from the book \cite{Samko1}:

$$
\Gamma(1 - \alpha) \ D^{\alpha}_x [F] = \frac{F(0)}{x^{\alpha}} + \int_0^x \frac{dF(t)}{(x-t)^{\alpha}} =
\int_0^x \frac{dF(t)}{(x-t)^{\alpha}},
$$
 since $ F(0) = F(0+) = 0. $ Therefore

$$
S_2  = {\bf E} (x - \xi)^{-2 \alpha} \ I(\xi < x) = \int_0^x \frac{d F(t)}{(x-t)^{2 \alpha}} = \Gamma(1 - 2 \alpha) \
D^{2 \alpha}[F](x).   \eqno(2.11)
$$
 It remains to substitute into equality (2.9),  taking into account the relation (1.4a). \par
 The asymptotical normality our estimate follows now from  the  classical one-dimensional CLT. \par

\vspace{3mm}

{\bf Remark 2.1.} As follows from the relation (2.8), under condition $ |G^{(\alpha)}(x)| < \infty  $
the variable $ f_{\alpha,\xi}(x) $ has a finite absolute expectation:

$$
{\bf E} |f_{\alpha,\xi} (x)| < \infty.
$$
 Therefore, on the basis of the Law of Large Numbers, the estimate $  G_{\alpha,n}(x) $ is consistent
with probability one only under the condition $ |G^{(\alpha)}(x)| < \infty.  $  \par

\vspace{3mm}

{\bf Remark 2.2.} We are not sure that offered in this report estimate  $ G_{\alpha,n}(x) $ of the value
 $ G^{(\alpha)}(x) $ is optimal, in the contradiction to the Kolmogorov's estimate of the ordinary distribution
function. \par

\vspace{3mm}

{\bf Remark 2.3.} Emerging in the theorem 2.1 the variable $  G^{( 2 \alpha)}(x),  $  which may be used by the
practical application, may be consistent estimated  as follows:

$$
 G^{( 2 \alpha)}(x) \approx G_{2 \alpha,n}(x),
$$
as  long as $ \alpha < 1/2. $ \par

\vspace{3mm}

{\bf Remark 2.4.} Non-asymptotical approach. \par

 Let $  x  $ be fixed positive number; we consider a non-asymptotical deviation

$$
P_{\alpha,x}(y) \stackrel{def}{=} \sup_n {\bf P} \left( \sqrt{n} \cdot \left| G_{\alpha,n}(x) - G^{(\alpha)}(x) \right| > y \right),  \ y \ge 3.
 \eqno(2.12)
$$

 Note that the summand  r.v. $  f_{\alpha,\xi} = x^{-\alpha} - (x - \xi)^{-\alpha}\cdot I(x > \xi)  $ has a heavy  tail. Namely, if

$$
F(x) - F(0) = F(x) \sim C_1 \ x^{\Delta}, \ x \to 0+, \ \Delta = \const \in (0,1], \ \eqno(2.13)
$$
then

$$
{\bf P} \left( \left| G_{\alpha,1}(x) - G^{(\alpha)}(x) \right| > y \right) \sim C_2(\alpha,x) \cdot y^{-\Delta/\alpha}, \ y \to \infty. \eqno(2.14)
$$
 The non-asymptotical bounds for $ \sqrt{n}  $ normed deviations of sums of these variables are obtained, e.g. the articles \cite{Adler1},
 \cite{Bravernan1},  \cite{Bravernan2}, \cite{Ostrovsky3}. We deduce the upper bound for considered probability using these results:

$$
P_{\alpha,x}(y) \le  C_3(\alpha,x) \cdot y^{-\Delta/\alpha} \ \ln y, \ y > 3. \eqno(2.15)
$$
 The lower estimate for these probability is trivial: as $  y \to \infty $

$$
P_{\alpha,x}(y) \ge  {\bf P} \left( \left| G_{\alpha,1}(x) - G^{(\alpha)}(x) \right| > y \right) =
C_4(\alpha,x) \cdot y^{-\Delta/\alpha}(1 +  o(1)  ). \eqno(2.16)
$$
  The ultimate value of the degree of the value $ \ln y $ in (2.15) is now unknown. \par

\vspace{4mm}

\section{Main result: error estimation in Lebesgue-Riesz norm.}

 \vspace{4mm}

 As long as the function $ D^{\alpha}[G](\cdot) $ and correspondingly its estimate  $ G_{\alpha,n}(\cdot)  $ both are discontinuous
 and all the more so are unbounded,
 we can not do the error estimation in the uniform norm, and still  can not apply  the CLT  in the Prokhorov-Skorokhod space,
 in contradiction  to the classical Kolmogorov's theorem.\par

 \vspace{3mm}

 We intent to investigate the $  L_q(dF) $ deviation of empirical derivative reliability function $ G_{\alpha,n} $
from its true value $ D^{\alpha}G:  \ W_{q,n} :=$

$$
W_{q,n} [ D^{\alpha}[G](\cdot), G_{\alpha,n}(\cdot) ] \stackrel{def}{=} \sqrt{n} \cdot
\left[ {\bf E} \int_R | D^{\alpha}[G](x) - G_{\alpha,n}(x)|^q \ dF(x) \right]^{1/q}, \eqno(3.1)
$$

 As usually, in order to evaluate the variable $  L_q, $ we need to establish the Central Limit Theorem in the Lebesgue-Riesz
space $  L_q(dF). $  \par

 Note first of all that the expression  for $ W_q $ in (3.1) does not depend on the function $  F(\cdot) $ inside the set
$ K = \{  F \}, $ as in the Kolmogorov's theorem; therefore we can and will suppose $ F(x) = x, 0 \le x \le 1, $  i.e.
$  a = 0, b = 1 $ and the r.v. $  \{  \xi_i \} $ have the uniform distribution on the set $  [0,1].$ \par

\vspace{3mm}

 We introduce some new notations.  $ x \cap y := \min(x,y), \ 1 \le q = \const < 1/\alpha,  $

$$
R_{\alpha}(x,y) \stackrel{def}{=} \frac{( x \cap y )^{1 - 2 \alpha}}{1 - 2 \alpha} \  - \frac{(x \ y)^{1 - \alpha}}{(1 - \alpha)^2}, \
0 \le x,y < 1, \ 0 < \alpha < 1/2, \eqno(3.2)
$$

$$
\sigma^2_{\alpha}(x):= R_{\alpha}(x,x) = \frac{ x^{1 - 2 \alpha}}{1 - 2 \alpha} \  - \frac{x^{2 - 2\alpha}}{(1 - \alpha)^2},
$$

$$
g_{\alpha}(x) = x^{-\alpha} - \frac{x^{1 - \alpha}}{1 - \alpha},
$$

\vspace{3mm}

$$
\Gamma(1 - \alpha) \ \zeta_{n}^{(\alpha)}(x) =  \Gamma(1 - \alpha)  \  \zeta_n(x) :=
$$

$$
n^{-1/2} \sum_{i=1}^n \left\{ f_{\alpha,\xi_i}(x) - \left[ x^{-\alpha} - \frac{x^{1 - \alpha}}{1 - \alpha}  \right]    \right\}  =
$$

$$
n^{-1/2} \sum_{i=1}^n \left\{ f_{\alpha,\xi_i}(x) -  g_{\alpha}(x) \right\},  \eqno(3.3)
$$

$ \zeta_{\infty}^{(\alpha)}(x) = \zeta^{(\alpha)}(x) = \zeta(x), \ 0 \le x \le 1 $ be a separable (moreover, continuous with probability one)
centered Gaussian random process with covariation function $ R_{\alpha}(x,y): $ \par

$$
\Cov(\zeta_{\alpha}(x), \zeta_{\alpha}(y)) = {\bf E} \zeta_{\alpha}(x) \cdot \zeta_{\alpha}(y)= R_{\alpha}(x,y). \eqno(3.4)
$$

 The ordinary  Lebesgue-Riesz space $ L_q(dF) = L_q(R_+,dF)  $ consists by definition on all the measurable functions $ f: R_+ \to R $ with
finite norms

$$
||f||_q =  ||f||L_q(dF) \stackrel{def}{=} \left[ \int_0^{\infty} |f(x)|^q \ dF(x) \right]^{1/q}, \ q = \const \ge 1.
$$

\vspace{4mm}

{\bf Theorem 3.1.}  Let $  F(\cdot) \in K $ and let $ 1 \le  q < 1/\alpha.  $ Our statement: the sequence of distributions generated in the space
$ L_q(dF)  $  by the random processes  $ \zeta_{n}^{(\alpha)}(\cdot) $ converges weakly as $ n \to \infty  $ to the random process
$ \zeta_{\infty}^{(\alpha)}(x) $ (the CLT in the space $ L_q(R_+, dF). $

\vspace{3mm}

{\bf Proof.} Note first of all that here

$$
G^{(\alpha)}(x) = D^{\alpha}[G](x) =
\frac{1}{\Gamma(1 - \alpha)} \cdot \left[ x^{-\alpha} - \frac{x^{1 - \alpha}}{1 - \alpha}  \right]. \eqno(3.5)
$$
as long as in this section $ G(x) = 1 - x, \ x \in (0,1). $  Therefore, all the processes $ \zeta_n(x) $ are centered.\par
 As before, it is sufficient to consider the centered random process $ \zeta_{1}^{(\alpha)}(x), \ x \in (0,1).  $ It is easy
to calculate its covariation function; it coincides with $  R_{\alpha}(x,y). $  \par
 It remains to establish the CLT in the Lebesgue-Riesz space $ L_q(0,1) $ for the sequence  $ \zeta_{n}^{(\alpha)}(\cdot). $
The using for us version of CLT in this  spaces is obtained, for example, in the fundamental monograph \cite{Ledoux1}, pp. 308-319.
Namely, the sufficient condition

$$
{\bf E} ||\zeta_1||_q^q < \infty \eqno(3.6)
$$
is here satisfied, as long as $ q < 1/\alpha. $ \par

 In  detail, let us denote

$$
K(\alpha,q) := \frac{2^{1 - 1/q}}{\Gamma(1 - \alpha)} \cdot \left[(1 - \alpha)^{-q} + (1 - \alpha q)^{-1} \right]^{1/q} < \infty,
\eqno(3.7)
$$
then

 $$
 \Gamma(1 - \alpha) \ |\zeta_1(x)| \le \frac{x^{1 - \alpha}}{1 - \alpha} +  (| (x - \xi)|I(x> \xi)|)^{-\alpha};
 $$

$$
\Gamma^q(1 - \alpha) \ |\zeta_1(x)|^q \le 2^{q-1} \left[ (1 - \alpha)^{-q} + x^{ - \alpha q } \right];
$$

$$
\Gamma^q(1 - \alpha) \ ||\zeta_1(\cdot)||_q^q = \Gamma^q(1 - \alpha) \ \int_0^1|\zeta_1(x)|^q \ dx
 \le 2^{q-1} \left[ (1 - \alpha)^{-q} +  (1 - \alpha q)^{-1} \right];
$$
so

$$
||\zeta_1(\cdot)||_q \le
\frac{2^{1 - 1/q}}{\Gamma(1 - \alpha)} \cdot \left[(1 - \alpha)^{-q} + (1 - \alpha q)^{-1} \right]^{1/q}
= K(\alpha,q) < \infty, \eqno(3.7a)
$$
since $ 0 < \alpha < 1, \ 1 \le q < 1/\alpha. $ \par

\vspace{3mm}

 This completes the proof of theorem 3.1. \par

\vspace{3mm}

{\bf Remark 3.1.} Note that the obtained estimate (3.7a) is deterministic, i.e.  is true still without the expectation {\bf E.} \par

\vspace{3mm}

{\bf Remark 3.2.} Let us denote

$$
Q_{\alpha}^{(q)}(u) = {\bf P} \left( || \zeta_{\alpha}(\cdot)||L_q(0,1) > u  \right), \ u = \const > 0. \eqno(3.8)
$$
 We deduce as a consequence of the theorem 3.1 for the values  $  q = \const \in [1, 1/\alpha),\ \alpha \in (0,1/2) $
and $  u > 0 $

$$
\lim_{n \to \infty} {\bf P} \left( || G_{\alpha,n}(\cdot)  - G^{(\alpha)}||L_q(dF) > \frac{u}{\sqrt{n} \Gamma(1 - \alpha)} \right) =
Q_{\alpha}^{(q)}(u), \eqno(3.9)
$$
therefore for sufficiently greatest values $  n  $

$$
 {\bf P} \left( || G_{\alpha,n}(\cdot)  - G^{(\alpha)}||L_q(dF) > \frac{u}{\sqrt{n} \Gamma(1 - \alpha)} \right) \approx
Q_{\alpha}^{(q)}(u). \eqno(3.10)
$$

 The asymptotical behavior of the probability $  Q_{\alpha}(u) $ as $ u \to \infty  $ is known, see
\cite{Piterbarg1},\cite{Piterbarg2}. Briefly, let us denote $ q' = q/(q - 1), \ q > 1, $ and introduce
the variable

$$
\beta := \sup_{h: ||h||L(q') = 1}  \left\{ \int_0^1   \int_0^1 R_{\alpha}(x,y) h(x) h(y) dx dy \right\},
$$
which may be computed in turn through solving of some non-linear integral equation; then

$$
\ln Q_{\alpha}^{(q)}(u)  \sim - \frac{u^2}{2 \beta^2},  \ u \to \infty.
$$
 The non-asymptotical estimates of this probability is obtained in \cite{Ostrovsky1}, chapter 4,
section 4.8; see also  \cite{Ostrovsky2}, chapter 3.\par

 It is clear that the equality (3.10) may be used by construction of confidence region for the unknown
function $  G^{(\alpha)}(\cdot) $  in the Lebesgue-Riesz norm $ L_q(dF) $ and for the testing of non-parametrical
hypotheses. \par

\vspace{3mm}

{\bf Remark 3.3.  Verification.} It is interest to note that on the case $ \alpha = 0, $ more exactly when $ \alpha \to 0+, $
the obtained before results coincide with the classical belonging to Kolmogorov, Mises  etc.

\vspace{4mm}

\section{Non-asymptotical error estimation in the Lebesgue-Riesz norm.}

\vspace{4mm}

 We intent to obtain in this section the non-asymptotical {\it upper} estimate for the supremum of loss function

$$
\overline{W}_q \stackrel{def}{=} \sup_n W_{q,n} [ D^{\alpha}[G](\cdot), G_{\alpha,n}(\cdot) ] \eqno(4.0)
$$
and as a consequence by means of Tchebychev's inequality  the  probability

$$
\sup_n  {\bf P} \left( || G_{\alpha,n}(\cdot)  - G^{(\alpha)}||L_q(dF) > \frac{u}{\sqrt{n} \Gamma(1 - \alpha)} \right) =:
 \overline{Q}_{\alpha}(u). \eqno(4.1)
$$

 We retain in this section all the notations and restrictions of third section; for instance, $ F(\cdot) \in K;  $ therefore
we can and will suppose that the r.v. $ \{  \xi(i) \} $  are independent and uniformly distributed on the set $ [0,1]. $\par

 The case $ \alpha = 0, $ i.e. when we consider the classical problem of estimation of ordinary distribution function  $ F(x) $
or equally the reliability function  $ G(x) = 1 - F(x) $
by means of empirical distribution function $ F_n(x), $ is investigated, and at once in the multidimensional case,
even in the uniform norm, i.e. when formally $ q = \infty,  $   in the work of
J.Kiefer \cite{Kiefer1}; more exact estimate see in the article \cite{Gaivoronsky1}. Indeed,

$$
{\bf P} \left(  \sqrt{n} \sup_{x} |G_n(x) - G(x)|> u  \right) \le 2 e^{- 2 u^2}, \ u \ge 1, \eqno(4.2)
$$
i.e. the exponential bound for $ \sqrt{n} $ normed uniform deviation $  \sup_{x} |G_n(x) - G(x)|.   $ \par
 Note first of all that the $ \sqrt{n}  $ exponential tail distributed confidence  region  for $ G^{(\alpha)}(\cdot)  $
based on our estimate $ G_{\alpha,n}(\cdot)  $ in the $ L_p  $ norm  is impossible when $ p \ge 1/\alpha,  $ in contradiction
to the classical ordinary case $  \alpha = 0. $ Namely,
we  can deduce the following  simple  {\it lower  } bound for $ \overline{W}_q  $

$$
\overline{W}_q \ge W_{q,1} [ D^{\alpha}[G](\cdot), G_{\alpha,1}(\cdot) ],
$$
and it is easily to calculate analogously to the relations (3.7) - (3.7a) that

$$
\Gamma(1 - \alpha) \cdot {\bf E}|| \left[ D^{\alpha}[G](\cdot) - G_{\alpha,n}(\cdot) \right] ||_q \ge \frac{C_1(\alpha)}{(1 - \alpha q)^{1/q} },
 \ q < 1/\alpha, \eqno(4.3)
$$
and

$$
\Gamma(1 - \alpha) \cdot {\bf E}|| \left[ D^{\alpha}[G](\cdot) - G_{\alpha,n}(\cdot) \right] ||_q  = \infty, \ q \ge 1/\alpha, \eqno(4.4)
$$

\vspace{3mm}

 We are going now to the obtaining of {\it upper} estimates for the value  $ \overline{W}_q. $ We suppose in the sequel
$ 2 \le q  < 1/\alpha. $\par

\vspace{3mm}

{\bf Theorem 4.1.} We propose under formulated above conditions $ 0 < \alpha < 1/2, \ 2 \le q < 1/\alpha $ etc.

$$
\Gamma(1 - \alpha) \cdot \sup_n W_{q,n} [ D^{\alpha}[G](\cdot), G_{\alpha,n}(\cdot) ] \le \frac{C_2(\alpha)}{(1 - \alpha q)^{1/q}}, \eqno(4.5)
$$
where $  C_2(\alpha) $ is continuous positive function on the closed segment $  \alpha \in [0, 1/2]. $ \par

\vspace{3mm}

{\bf Proof.} Denote

$$
\tau_i(x) = f_{\alpha,\xi(i)}(x) - {\bf E} f_{\alpha,\xi(i)}(x), \eqno(4.6)
$$

$$
S_n(x) = \frac{1}{\sqrt{n}} \sum_{i=1}^n \tau_i(x), \eqno(4.7)
$$
then $ \tau_i(x)  $ is a sequence of independent identical distributed centered random fields which are proportional with
coefficient $ 1/\Gamma(1 - \alpha) $ to the considered before r.f. $  \zeta_i(x). $ \par

 Let us consider the sequence of random variables

$$
V_{\alpha,q}(n) = ||S_n(\cdot)||_q^q  = \int_0^1 |S_n(x)|^q \ dx,
$$
then we have using Fubini-Tonelli theorem under our condition  $ 2 \le q < 1/\alpha $

$$
{\bf E} V_{\alpha,q} (n) = \int_0^1 {\bf E} |S_n(x)|^q \ dx=
\int_0^1 {\bf E} \left| n^{-1/2} \sum_{i=1}^n \tau_i(x)  \right|^q \ dx. \eqno(4.8)
$$

We intent to exploit the famous Rosenthal's inequality, see \cite{Rosenthal1}, \cite{Ostrovsky4}.
Namely,  for  arbitrary sequence $  \{ \zeta_k \}  $ of  independent centered random variables

$$
 \left|\ n^{-1/2} \ \sum_{k=1}^n \zeta_k \ \right|_q \le K_{R} \cdot \frac{q}{\ln q}  \cdot \sqrt{ \sum_{k=1}^n |\zeta_k|_q^2/n  },
 \hspace{4mm} q \ge 2, \eqno(4.9)
$$
where  the "Rosenthal's"  constant  $  K_{R} $  is less than 0.6535, see  \cite{Ostrovsky4}. \par

 If the r.v. $ \{  \zeta_k \}  $  are in addition identically distributed, then

$$
 \left|\ n^{-1/2} \ \sum_{k=1}^n \zeta_k \ \right|_q \le K_{R} \cdot \frac{q}{\ln q}  \cdot |\zeta_1|_q.
 \eqno(4.10)
$$

 As long as in this section $ 2 \le q < 1/\alpha, $

$$
\frac{q}{\ln q} \le \max \left[ \frac{2}{\ln 2}, \    \frac{1/\alpha}{|\ln \alpha|}  \right],
$$
and if we denote

$$
K_{R,\alpha} = K_R \cdot \max \left[ \frac{2}{\ln 2}, \    \frac{1/\alpha}{|\ln \alpha|}  \right], \eqno(4.11)
$$
then there holds the following inequality (under our conditions)

$$
 \left|\ n^{-1/2} \ \sum_{k=1}^n \zeta_k \ \right|_q \le K_{R,\alpha} \cdot |\zeta_1|_q. \eqno(4.12)
$$

 We  conclude after substitution into (4.8)

$$
{\bf E} V_{\alpha,q} (n) \le K(\alpha,q)^q \cdot  K^q_{R,\alpha}.  \eqno(4.13)
$$

 It remains to extract the root of degree $  q  $ from last inequality to obtain the estimate

$$
\sup_n W_{q,n} [ D^{\alpha}[G](\cdot), G_{\alpha,n}(\cdot) ] \le K(\alpha,q) \cdot  K_{R,\alpha},
$$
which is equivalent to the assertion if theorem 4.1 with explicit evaluate of constant. \par

\vspace{4mm}

\section{Estimation of fractional derivatives of spectral function.}

 \vspace{4mm}

 Let us consider in this section the classical problem of spectral density estimation. \par

  Let $ \eta_k, \ k = 1,2,\ldots,n $ be real valued, centered: $  {\bf E} \eta(k) = 0 $ Gaussian distributed stationary random
 sequence (process) with (unknown) even covariation function $  r = r(m), $  spectral function
 $ F(\lambda), \ \lambda \in [0, 2 \pi],  \ F(0+) = F(0) = 0, $ and with  spectral density $  f(\lambda), $ (if there exists):

$$
r(m) = \Cov(\eta(j+m), \eta(j)) =  {\bf E} \eta(j+m)\cdot \eta(j) =
$$

$$
\int_{[0, 2 \pi]} \cos(\lambda m) dF(\lambda) =
\int_{[0, 2 \pi]} \cos(\lambda m) f(\lambda) \ d \lambda, \eqno(5.1)
$$
so that

$$
F(\lambda) = \int_0^{\lambda} f(t) dt = I^{(1)}[f](\lambda).
$$

 The periodogram of this sequence will be denoted by $ J_n(\lambda), \ 0 \le \lambda \le 2 \pi: $

$$
J_n(\lambda) := (2 \pi n)^{-1} \left| \sum_{k=1}^n e^{i k \lambda} \eta(k)  \right|^2; \ i^2 = - 1. \eqno(5.2)
$$

 We intent here to estimate the fractional derivative  $  F^{(\alpha)}(\lambda)  $ of the spectral function  $  F(\lambda). $ \par
 Recall that the problem of $  F(\cdot) $ estimation is well posed, in contradiction to the problem of spectral density
 $ f(\cdot) = F^{(1)}(\cdot) $  estimation.\par
 We assume as before $ 0 < \alpha < 1/2,  $ and denote $  \beta = 1 - \alpha; \ \beta \in (1/2, 1). $  \par

Heuristic arguments. We have using the group properties of the fractional derivative-integral operators

$$
F^{(\alpha)} = D^{\alpha}[F] = D^{\alpha} I^1[f] =D^{\alpha} D^{-1}[f] = D^{\alpha - 1}[f] =
$$

$$
I^{1 - \alpha}[f] = I^{\beta}[f] \approx  I^{\beta} [J_n].\eqno(5.3)
$$

 Thus, we can offer as an estimation of $ F^{(\alpha)} $ the following statistics

$$
F_{\alpha,n}(\lambda) :=  I^{\beta} [J_n](\lambda) = \frac{1}{\Gamma(\beta)}\int_0^{\lambda} \frac{J_n(t) \ dt}{(t - \lambda)^{1-\beta }} =
\frac{1}{\Gamma(1-\alpha)}\int_0^{\lambda} \frac{J_n(t) \ dt}{(t - \lambda)^{\alpha }}. \eqno(5.4)
$$

\vspace{3mm}

{\bf Theorem 5.1.} Suppose as before $ 0 < \alpha < 1/2 $ and
that the spectral density  $  f(\lambda) $ there exists and is continuous on the circle $ [0, 2 \pi], $
in particular

$$
f(0) = f(0+) = f(2 \pi -0)  = f(2 \pi).
$$

  Then the estimate $ F_{\alpha,n}(\lambda)   $  of the fractional derivative $ F^{(\alpha)}(\lambda)  $ is asymptotically
unbiased:

$$
{\bf E} F_{\alpha,n}(\lambda) =  F^{(\alpha)}(\lambda) + O(1/n), \eqno(5.5)
$$
and

$$
\lim_{n \to \infty} n \cdot \Var \left[ F_{\alpha,n}(\lambda)\right] = \frac{4 \pi}{\Gamma^2(1 - \alpha)} \cdot \int_0^{\lambda}
\frac{f^2(\nu) d \nu}{ (\lambda - \nu)^{2 \alpha} } =
$$

$$
\frac{4 \pi \Gamma(1 - 2 \alpha)}{\Gamma^2(1 - \alpha)} \cdot I^{2 \alpha}[f^2](\lambda) < \infty. \eqno(5.6)
$$

 Note that the last integral is finite since the function $ f $ is bounded and  $  \alpha < 1/2. $ \par

 More generally,

$$
\lim_{n \to \infty} n \cdot \Cov \left\{ F_{\alpha,n}(\lambda), F_{\alpha,n}(\mu)  \right\} =
\frac{4 \pi \Gamma(1 - 2 \alpha)}{\Gamma^2(1 - \alpha)} \cdot D^{-\alpha}_{\lambda} \left[ D^{-\alpha}_{\mu}[f^2]  \right]=
$$

$$
\frac{4 \pi}{\Gamma^2(1 - \alpha)} \cdot
\int_0^{ \lambda \cap \mu } \frac{f^2(\nu) \ d \nu}{(\lambda - \nu)^{\alpha} \ (\mu - \nu)^{\alpha}}
=: \Theta_{\alpha}(\lambda, \mu). \eqno(5.7)
$$

\vspace{3mm}

{\bf Proof.} Our  assertion follows immediately from
 the following proposition, see  the  fundamental monograph of T.W.Anderson
\cite{Anderson1}, chapter 5, page 564-572, theorem 9.3.1:
if $ w(\lambda, \nu) $ is non-negative integrable function, then

$$
\lim_{n \to \infty} \int_0^{2 \pi} w(\lambda, \nu) J_n(\nu) d \nu = \int_0^{2 \pi} w(\lambda, \nu) f(\nu) d \nu, \eqno(5.8a)
$$

$$
\lim_{n \to \infty} n \cdot \Cov \left( \int_0^{2 \pi} w(\lambda_1, \nu) J_n(\nu) d \nu,
\int_0^{2 \pi} w(\lambda_2, \nu) J_n(\nu) d \nu  \right) =
$$

$$
4 \pi \int_0^{2 \pi} w(\lambda_1, \nu) \ w(\lambda_2, \nu) \ f^2(\lambda) \ d \lambda, \eqno(5.8b)
$$
with remainder terms. We choose $ w(\lambda, \nu) = |\lambda - \nu|^{-\alpha}; $ it is easy to verify that
all the conditions of the mentioned result are satisfied. \par
 Recall also that the considered stationary sequence $  \{ \eta(k) \}  $ is Gaussian, i.e. without cumulant function. \par

\vspace{3mm}

{\bf Remark 5.1.} Emerging in the equality (5.6) the variable

$$
 I^{2 \alpha}[f^2](\lambda) = \frac{1}{\Gamma(2 \alpha)} \int_0^{\lambda} \frac{f^2(\nu) d \nu}{(\lambda - \nu)^{1 - 2 \alpha}}
$$
may be $ n^{-1/2} \ - $ consistent estimated as follows:

$$
 I^{2 \alpha}[f^2](\lambda) \approx \frac{1}{\Gamma(2 \alpha)} \int_0^{\lambda} \frac{J_n^2(\nu) d \nu}{(\lambda - \nu)^{1 - 2 \alpha}}. \eqno(5.9)
$$

\vspace{3mm}

{\bf Remark 5.2.} I.A.Ibragimov in  \cite{Ibragimov1} proved the asymptotical normality of the random process
$ \sqrt{n} \left\{I^1[J_n](\lambda) - I^1[F](\lambda) \right\}  $ as $  n \to \infty $ in the space $ C(0, 2 \pi) $ of continuous
functions.  See also \cite{Dahlhaus1}, \cite{Levit1}. A fortiori, the sequence of random processes

$$
\zeta_n(\lambda) =  \sqrt{n} \cdot \left\{ F_{\alpha,n}(\lambda) - F^{(\alpha)}(\lambda)  \right\}
$$
converges weakly in the space $   C(0, 2 \pi)  $ as $ n \to \infty $ to the  centered separable Gaussian process
$  \zeta_{\infty} $ with covariation function $  \Theta_{\alpha}(\lambda, \mu). $ Therefore

$$
{\bf P} \left(  \sqrt{n} \cdot \max_{\lambda} \left| \left\{ F_{\alpha,n}(\lambda) - F^{(\alpha)}(\lambda) \right\} \right|  > u  \right)
\approx {\bf P} (\max_{\lambda} | \zeta_{\infty}(\lambda) | > u), \ u = \const > 0. \eqno(5.10)
$$

 The asymptotical as $  u \to \infty  $ behavior of the last probability is fundamental investigated in the monograph
\cite{Piterbarg2}, see also  \cite{Piterbarg1}:

$$
 {\bf P} (\max_{\lambda} | \zeta_{\infty}(\lambda) | > u) \sim H(\alpha) \ u^{\kappa - 1} \ \exp \left(-u^2/\sigma^2 \right),\eqno(5.11)
$$

$$
H(\alpha), \ \kappa = \const, \ \sigma^2 = \sigma^2(\alpha) = \max_{\lambda \in (0, 2 \pi)} \Theta_{\alpha}(\lambda,\lambda). \eqno(5.11a)
$$

 The last equalities may be used by construction of confidence region  for $ F^{(\alpha)}(\cdot) $ in the uniform norm.
Indeed, let $ 1 - \delta $  be the reliability of confidence region, for example, $ 0.95  $ or $ 0.99 $ etc.
Let $ u_0 = u_0(\delta) $ be a maximal root of the equation

$$
H(\alpha) \ u_0^{\kappa - 1} \ \exp \left(-u_0^2/\sigma^2 \right) = \delta,
$$
then with probability $ \approx 1 - \delta $

$$
 \sup_{\lambda \in (0, 2 \pi)} \left| F_{\alpha,n}(\lambda) - F^{(\alpha)}(\lambda)  \right| \le \frac{u_0(\delta)}{\sqrt{n}}. \eqno(5.12)
$$

\vspace{4mm}

\section{Multidimensional case.}

\vspace{3mm}

 We consider in this section the problem of statistical estimates of fractional derivative for multidimensional distribution
function.  We restrict ourselves for simplicity only two-dimensional case $  d = 2. $ \par
 In detail, let $ \{ (\xi(i), \eta(i)) \}, i = 1,2,\ldots,n $ be a two dimensional non - negative sample with common distribution
function. We define the reliability function $ G = G(x,y) $  as follows:

$$
G(x,y) = {\bf P}(\xi(i) \ge x, \ \eta(i) > y), \ x,y \ge 0.\eqno(6.1)
$$
 Let $  \alpha, \beta  = \const $ be two numbers such that $  0 < \alpha, \beta < 1; $ (we will suppose further that
$  0 < \alpha, \beta < 1/2.) $ The partial mixed fractional derivative $ D^{ \alpha, \beta}_{x,y}[G](x,y)   $
again of  Rieman-Liouville type
of order $ (\alpha, \beta) $ of a function $ G(\cdot, \cdot) $ at the positive points $ (x,y)  $ is defined as  follows:

$$
G^{(\alpha, \beta)}(x,y) =  D^{ \alpha, \beta}_{x,y}[G](x,y)  \stackrel{def}{=}   D^{\alpha}_x D^{\beta}_y [G] =
\frac{1}{\Gamma(1 - \alpha)} \frac{1}{\Gamma(1 - \beta)}\times
$$

$$
\frac{\partial^2}{ \partial x \partial y } \int_0^x \int_0^y \frac{G(t,s)\ dt \ ds}{ (x-t)^{\alpha} (y-s)^{\beta} }, \eqno(6.2)
$$
see, e.g. \cite{Samko1}, chapter 24. We put as before $ D^{ \alpha, \beta}_{x,y}[G](x,y) = 0 $ if at the point (in the plane) $ (x,y) $
the expression (4.2) for $ D^{ \alpha, \beta}_{x,y}[G](x,y) $ does not exists. \par

 Note that in general case $  D^{\alpha}_x D^{\beta}_y [H] \ne D^{\beta}_y D^{\alpha}_x [H],  $ but if the function $ G = G(x,y) $
is factorable: $ H(x,y) = g_1(x) g_2(y)  $ and both the functions  $ g_x(\cdot)  $ and $ g_2(\cdot) $ are "differentiable" at the
points $  x  $ and $  y  $ correspondingly:

$$
\exists D^{\alpha}[g_1](x), \hspace{6mm} \exists D^{\beta}[g_2](y),
$$
then really

$$
D^{\alpha}_x D^{\beta}_y [H] = D^{\beta}_y D^{\alpha}_x [H] = D^{\alpha}_x[g_1](x) \cdot D^{\beta}_y [g_2](y). \eqno(6.3)
$$

 Introduce as a capacity of the function  $  H  $ the following:

$$
H(x,y) := I(\xi \ge x, \ \eta \ge y), \ \xi = \xi(1), \eta = \eta(1).
$$
 As long as

$$
 I(\xi \ge x, \ \eta \ge y) = I(\xi \ge x) \cdot I(\eta \ge y),
$$
the function $  H = H(x,y) $ is factorable and therefore (see (6.3))

$$
\Gamma(1 - \alpha) \Gamma(1 - \beta) D^{\alpha}_x D^{\beta}_y [H] = f_{\alpha,\xi}(x) \cdot f_{\beta,\eta}(y). \eqno(6.4)
$$

\vspace{3mm}

 The consistent with probability  one in each fixed point $  (x,y) $ estimate of the function $ G^{(\alpha, \beta)}(x,y) $
is follows:

$$
\Gamma(1 - \alpha) \Gamma(1 - \beta)   G_{\alpha, \beta,n }(x,y) \stackrel{def}{=} n^{-1} \sum_{i=1}^n f_{\alpha,\xi(i)}(x) \cdot f_{\beta,\eta(i)}(y).
\eqno(6.5)
$$

 It is easily to verify  that the estimate  $ G_{\alpha, \beta,n }(x,y) $ obeys at the same properties as its one-dimensional predecessor
$ G_{\alpha, n }(x),  $  for example, is unbiased, satisfies LLN and CLT. \par

 Note that despite the function

$$
h_{\alpha, \beta}(x,y) = h_{\alpha, \beta}(\xi, \eta; x,y) :=  f_{\alpha,\xi}(x) \cdot f_{\beta,\eta}(y), \ \xi = \xi(1), \eta = \eta(1)\eqno(6.6)
$$
is also factorable, we do not suppose the independence  of the r.v. $  (\xi, \eta). $\par

 Define as before the following sequence of random fields

$$
S_n(x,y) = \frac{1}{\sqrt{n}} \sum_{i=1}^n \left( h_{\alpha, \beta}(\xi(i), \eta(i); x,y) -
\Gamma(1 - \alpha) \Gamma(1 - \beta) G^{(\alpha, \beta)}(x,y)  \right), \eqno(6.7)
$$
so that

$$
\Gamma(1 - \alpha) \Gamma(1 - \beta) \left[   G_{(\alpha, \beta,n)}(x,y)  -  G^{(\alpha, \beta)}(x,y) \right] = n^{-1/2} S_n(x,y). \eqno(6.8)
$$

\vspace{3mm}

 Further, it is clear that

$$
|| h_{\alpha, \beta}(\cdot, \cdot)||_{q,X; r,Y} \asymp (1 - \alpha q)^{-1/q} \ (1 - \beta r)^{-1/r}, \ 1 \le q < 1/\alpha,
1 \le r < 1/\beta. \eqno(6.9)
$$

\vspace{3mm}

 Let return to the source problem and let us consider only the non-mixed case $ r = q; $ then
 $$
 ||h(\cdot, \cdot)||_{q,q} =  ||h(\cdot, \cdot)||_{q, F} =   \left[\int_X \int_Y  |h(x,y)|^q  F(dx,dy) \right]^{1/q}.\eqno(6.10)
 $$
Assume also $  q \ge 2, $ so $ 2 \le q = r < \min(1/\alpha, 1/\beta). $\par

We will distinguish two essentially different variants: $ V_1: \beta < \alpha $ and $ V_2: \beta = \alpha. $
 The case $  \beta > \alpha  $ may be considered analogously. \par

\vspace{3mm}

 {\bf  First possibility:} $  \beta < \alpha. $ \par

\vspace{3mm}

 We find after simple  calculations  as in fourth section  using Rosenthal's inequality:

$$
 \sup_n {\bf E} || S_n(\cdot, \cdot)   ||_{q,F}^q \asymp C_1(\alpha, \beta, q) \cdot (1 - \alpha q)^{-1}, \ 1 \le q < 1/\alpha, \eqno(6.11)
$$
and

$$
 \sup_n {\bf E} || S_n(\cdot, \cdot)   ||_{q,F}^q  = \infty, \  q \ge 1/\alpha. \eqno(6.11a)
$$

\vspace{3mm}

 {\bf  Second possibility:} $  \beta = \alpha. $ \par

\vspace{3mm}

 We have analogously

$$
 \sup_n {\bf E} || S_n(\cdot, \cdot)   ||_{q,F}^q \asymp C_2(\alpha, q) \cdot (1 - \alpha q)^{-2}, \ 1 \le q < 1/\alpha, \eqno(6.12)
$$
and

$$
 \sup_n {\bf E} || S_n(\cdot, \cdot)   ||_{q,F}^q  = \infty, \ q \ge 1/\alpha. \eqno(6.12a)
$$

\vspace{4mm}

 As a consequence: the sequence of r.f. $ S_n(\cdot, \cdot)  $ in both the considered cases
 satisfies the CLT in the space $ L_q(F(dx, dy))  $  iff $  1 \le q < 1/\alpha. $ \par

\vspace{4mm}

 Thus, there is a possibility to built  the asymptotical and non-asymptotical confidence region for estimated
mixed fractional derivative $ F^{(\alpha, \beta)}(\cdot, \cdot), \ \alpha, \beta < 1/2 $  still in the multivariate case
 in the $ L_q(dF) $ norm as well as  in the fixed point $ (x_0, y_0). $ \par

\vspace{4mm}

\section{Estimation of fractional derivative in Grand Lebesgue Space norm.}

\vspace{3mm}

 Let  $ (X, M, \mu )  $ be a probability space with non-trivial probability measure $  \mu, $ and let also $ \psi = \psi(q), \ 1 \le q < s,
 \ s = \const \in (1, \infty]  $  be continuous on the open interval $   ( 1, s )  $ bounded from below  function.  By definition, a
Grand Lebesgue Space (GLS) $ G\psi $ over our triplet $ (X, M, \mu )  $ consists  on all the measurable functions $  f: X \to R $  with
finite norm

$$
||f||G\psi \stackrel{def}{=} \sup_{q \in (1,s)} \left[  \frac{|f|_q}{\psi(q)}  \right]. \eqno(7.1)
$$
 Hereafter

 $$
 |f|_q  = \left[ \int_X |f(x)|^q \ \mu(dx)  \right]^{1/q}
 $$
and we will denote $ s = \supp \psi. $ \par
 The detail investigation  of these spaces see, e.g. in \cite{Fiorenza1},  \cite{Fiorenza2},
\cite{Iwaniec1}, \cite{Iwaniec2}, \cite{Kozachenko1},  \cite{Ostrovsky1}, \cite{Ostrovsky7}, \cite{Ostrovsky8}.\par

\vspace{3mm}

We choose supposing  without loss of generality  $ a = 0, \ b = 1,  $ so that $  \dim \xi = 1 $
and $ F(0+) = 0, \ F(1 - 0) = 1, \ F \in K, \ \alpha \in (0,1/2), $

$$
X = [0,1] \otimes \Omega,
$$
so that the measure $ \mu $ is direct product of ordinary Lebesgue measure $  dx $ and probability measure $  P: $

$$
\mu(A \otimes B) = \int_A dx \cdot {\bf P}(B), \ A \subset [0,1], \ B \subset \Omega, \eqno(7.2)
$$
where $ \Omega  = \{  \omega \} $ is source probability space, i.e. in which is defined  our sample $ \{  \xi(i) \}. $ \par

\vspace{4mm}

{\bf Proposition 7.1.}

$$
\sup_n  \mu \{ (x,\omega): |S_n(x)| > u   \} \le C_3(\alpha) \ u^{-\alpha} \ \ln u, \ u \ge e. \eqno(7.3)
$$

{\bf  Proof.} Put

$$
\psi_{\alpha}(q) = (1 - \alpha q)^{-1/q}, \ 1 \le q < 1/\alpha.
$$
 Note that

$$
\psi_{\alpha}(q) \asymp (1 - \alpha q)^{-\alpha}, \ 1 \le q < 1/\alpha.
$$

The assertion (4.5) may be rewritten as follows.

$$
\sup_n |   S_n(\cdot, \cdot) |_{q, \mu} \le C_4(\alpha) \ \psi_{\alpha}(q), \eqno(7.4)
$$
or equally on the language of the Grand Lebesgue Spaces

$$
\sup_n || S_n(\cdot, \cdot) ||G\psi_{\alpha} \le C_4(\alpha)  < \infty. \eqno(7.5)
$$
 The tail estimate (7.3) follows immediately from one of results of the article \cite{Ostrovsky7};
 see also \cite{Ostrovsky8}.\par

\vspace{3mm}

 {\bf Remark 7.1.}  It is not hard to generalize this result into the  multidimensional case  described in the $  6^{th} $
section. Namely, if in the notations of the $  6^{th} $ section $  0 < \beta < \alpha < 1/2,  $ then

$$
\sup_n  \mu \{ (x,y, \omega): |S_n(x,y)| > u   \} \le C_5(\alpha, \beta) \ u^{-\alpha} \ \ln u, \ u \ge e;  \eqno(7.6)
$$
if $  0 < \beta = \alpha < 1/2,  $ then

$$
\sup_n  \mu \{ (x,y, \omega): |S_n(x,y)| > u   \} \le C_6(\alpha) \ u^{-\alpha} \ \ln^2 u, \ u \ge e.  \eqno(7.7)
$$

\vspace{4mm}

\section{Concluding remarks.}

\vspace{3mm}

{\bf A. Weight case. } \\

 Perhaps, it is interest to  investigate the error of the approximation of a form

 $$
 V(x) \cdot ( W \cdot F)^{(\alpha)}(x) \approx V(x) \cdot ( W \cdot F)_{n,\alpha}(x),
 $$
or analogously

$$
 V(x) \cdot ( W * F)^{(\alpha)}(x) \approx V(x) \cdot ( W * F)_{n,\alpha}(x),
$$
or analogously

$$
 V(x) * ( W * F)^{(\alpha)}(x) \approx V(x) *( W * F)_{n,\alpha}(x),
$$
where $ V(x), W(x)  $ are two weight functions, for instance, $ V(x) = |x|^{\gamma}, \ W(x) = |x|^{\Delta},  \hspace{6mm}
\gamma, \Delta = \const.  $ \par

\vspace{3mm}

{\bf B. Semi-parametric case.} \\

 Let $  \xi(i), \ i = 1,2,\ldots,n  $ be a sample of a volume $  n  $ with parametric family of regular distribution of a form

$$
{\bf P} (\xi(i) < x) = F(x; \theta),
$$
 where $ \theta \in \Theta \subset R^k, \ k < \infty $ is  (multidimensional, in general case) unknown numerical parameter. \par

 Denote by $ \hat{\theta}_n  $  the maximum likelihood estimate of the parameter $ \theta $ builded on our sample. The asymptotical
tail behavior of distribution  for the following statistic

$$
\tau_n := \sqrt{n} \cdot || F(\cdot, \hat{\theta}_n ) - F_n(\cdot) ||L,
$$
where $  || \cdot||L  $ is some Banach functional norm in the space $ x \in R, $ is in detail investigated in \cite{Piterbarg1},
\cite{Piterbarg2}, chapter 5.\par
  By our opinion, it is interest to obtain also the asymptotical  tail behavior of the following statistic

$$
\tau_{n,\alpha} := \sqrt{n} \cdot || F^{(\alpha)}(\cdot, \hat{\theta}_n ) - F_{n,\alpha}(\cdot) ||L,
$$
$  \alpha = \const \in (0, 1/2). $\par

\vspace{3mm}

{\bf C. Applications (possible) in  statistics.}\\

 The asymptotical  tail behavior of the  statistic $ \tau_{n,\alpha} $ may be used perhaps  in turn in statistics,
for instance, for the verification of semi-parametrical hypotheses and detection of distortion times etc. \par

\end{document}